\documentclass[12pt]{article}
\usepackage{graphicx, amsmath}
\usepackage{amsfonts}
\usepackage{verbatim}
\usepackage{latexsym}
\makeatletter
\@addtoreset{equation}{section}

\begin{document}

\newcommand{\E}{\mathbb{E}}
\newcommand{\PP}{\mathbb{P}}
\newcommand{\RR}{\mathbb{R}}

\newtheorem{theorem}{Theorem}[section]
\newtheorem{lemma}{Lemma}[section]
\newtheorem{coro}{Corollary}[section]
\newtheorem{defn}{Definition}[section]
\newtheorem{assp}{Assumption}[section]
\newtheorem{expl}{Example}[section]
\newtheorem{prop}{Proposition}[section]
\newtheorem{remark}{Remark}[section]

\newcommand\tq{{\scriptstyle{3\over 4 }\scriptstyle}}
\newcommand\qua{{\scriptstyle{1\over 4 }\scriptstyle}}
\newcommand\hf{{\textstyle{1\over 2 }\displaystyle}}
\newcommand\hhf{{\scriptstyle{1\over 2 }\scriptstyle}}

\newcommand{\eproof}{\indent\vrule height6pt width4pt depth1pt\hfil\par\medbreak}

\def\a{\alpha} \def\g{\gamma}
\def\e{\varepsilon} \def\z{\zeta} \def\y{\eta} \def\o{\theta}
\def\vo{\vartheta} \def\k{\kappa} \def\l{\lambda} \def\m{\mu} \def\n{\nu}
\def\x{\xi}  \def\r{\rho} \def\s{\sigma}
\def\p{\phi} \def\f{\varphi}   \def\w{\omega}
\def\q{\surd} \def\i{\bot} \def\h{\forall} \def\j{\emptyset}
\def\R{\mathbb R} 
\def\be{\beta} \def\de{\delta} \def\up{\upsilon} \def\eq{\equiv}
\def\ve{\vee} \def\we{\wedge}
\def\dd{\delta}
\def\F{{\cal F}}
\def\T{\tau} \def\G{\Gamma}  \def\D{\Delta} \def\O{\Theta} \def\L{\Lambda}
\def\X{\Xi} \def\Si{\Sigma} \def\W{\Omega}
\def\M{\partial} \def\N{\nabla} \def\Ex{\exists} \def\K{\times}
\def\V{\bigvee} \def\U{\bigwedge}

\def\1{\oslash} \def\2{\oplus} \def\3{\otimes} \def\4{\ominus}
\def\5{\circ} \def\6{\odot} \def\7{\backslash} \def\8{\infty}
\def\9{\bigcap} \def\0{\bigcup} \def\+{\pm} \def\-{\mp}
\def\la{\langle} \def\ra{\rangle}

\def\tl{\tilde}
\def\trace{\hbox{\rm trace}}
\def\diag{\hbox{\rm diag}}
\def\for{\quad\hbox{for }}
\def\refer{\hangindent=0.3in\hangafter=1}

\newcommand\wD{\widehat{\D}}
\newcommand{\ka}{\kappa_{10}}

\title{
\bf {Stability of Tamed EM scheme of Neutral Stochastic Differential Delay Equations}
}
\author
{ {\bf Yanting Ji}$^{\tt a}$\thanks{Contact e-mail address: mathjyt@hotmail.com}\ \, {\bf Guowei Wan}$^{\tt b}$ \ \, {\bf Zheng Zhou}$^{\tt a}$ 
	\\ [1ex]$^{\tt a}$ Department of Finance\\Zhejiang Shuren University, China\\
$^{\tt b}$Institute of Urban and Rural circulation Economy\\ Zhejiang Shuren University, China\\}

\date{}

\maketitle
\begin{abstract}
	In this paper, exponential mean-square stability and almost sure stability of the tamed EM scheme to neutral stochastic differential delay equation are investigated. Surprisingly, the exponential mean-square stability can reproduce the almost sure stability without placing extra conditions.
\\

	\smallskip
	
	
\end{abstract}
\section{Introduction}
Many real world phenomena could be modeled by stochastic dynamics systems with delay. For example, in \cite{Mao1}, Mao, Yuan and Zou modified the classic Lotka-Volterra model by adding a delay argument, by which, an extended range of phenomena could be covered. In fact, the principle of causality, the future state of a system is not dependent of the past states and is solely determined by the present, failed to hold. In this paper, a more   general form of stochastic differential equations(SDEs) is investigated, the neutral stochastic differential delay equations (NSDDEs), which describes a wide variety of natural and artificial dynamical system, which depend not only on present states but also on past states as well as the derivative of the delay variables. 

Most SDEs are non-linear and the corresponding explicit solutions can hardly be found. The NSDDEs, as a general from of SDEs, also share this property. Therefore, the numerical scheme for NSDDEs plays a crucial role in studying these stochastic dynamical systems. In evaluation the quality of a numerical scheme, we need to consider two important aspects: convergence and stability.

The most widely applied numerical scheme is the Euler-Maruyama (EM) scheme, which is the natural extension of finite difference scheme for ordinary differential equations to the stochastic counterparts. For several decades, the numerical scheme for stochastic differential equations has drawn a lot of attention, many convergence results have been developed. In \cite{Klo}, Kloeden and Platen showed that, if both drift and diffusion coefficients are globally Lipschitz, then the classic EM scheme converges to the explicit solution of SDEs. In \cite{Baker}, Baker and Buckwar developed an explicit numerical scheme for stochastic differential delay equations under global Lipschitz condition and linear growth condition.

It is the first time that, Higham, Mao and Stuart \cite{Hig} establish strong convergence results under super-linear condition and the moment boundedness condition. In the same paper, an open question was asked whether the moment of the EM scheme is bounded within finite time if the coefficients of an SDE are not globally Lipschitz continuous. In \cite{Sab} Gy\"ongy and Sabanis showed that the EM scheme is convergent under local monotonicity condition for the stochastic differential delay equations. Jacob, Wang and Yuan \cite{Jacob} further investigated the convergence of EM scheme to stochastic differential delay equations with jumps under similar conditions.  Bao and his co-authors \cite{Bao1} extended the EM scheme to a more general case where the delay argument is a time segment rather than a fixed point. Hutzenthaler, Jentzen and Kloeden \cite{Hu} proved that once the global Lipschitz condition was dropped the EM scheme may explode in the finite time. To overcome this difficulty, Hutzenthaler, Jentzen and Kloeden \cite{Jen} developed a new type of numerical scheme, the so-called tamed scheme, which set a bound for the drift term so that the numerical solutions would not explode within the finite time. Sabanis \cite{Sab} proved convergence of tamed scheme to stochastic differential delay equations, while Ji and Yuan \cite{Ji2} further extended the convergence result of the tamed scheme to the NSDDEs.

With the convergence result in hand, the stability property, for both stability of explicit and numerical solutions to NSDDEs, have been developed during the past several years, for instance \cite{Bao}, \cite{Lan}, \cite{Mao},and \cite{Wu}, here we only mention a few. 

In \cite{Lan}, the exponential stability of both exact and $\theta$-EM scheme of NSDDEs with local Lipschitz condition has been investigated. However, the $\theta$-EM scheme is an implicit scheme rather than explicit scheme, which could add extra computational effort to solve such an implicit numerical system. In \cite{Zong}, the convergence and stability of a  modified tamed scheme for SDEs with non-Lipschitz continuous coefficients. In \cite{Zong1}, Zong and Wu investigated the $p$-th moment and almost sure exponential stability of the exact and numerical solutions of NSDDEs by virtue of the Lyapunov method. However, the coefficients are required to satisfied regular conditions.

To our best knowledge, there is little known result focus on the exponential stability of the tamed numerical solutions to this NSDDEs. Enlighten by \cite{Lan} and \cite{Zong1}, in this paper, we investigate the exponential mean-square stability of tamed numerical solutions to NSDDEs, whose coefficients are non-Lipschitz. Hence, we need to overcome more technical difficulties.

The rest of paper is organized as follows: Section $2$ gives some preliminary results. In Section $3,$ the exponential mean-square stability for tamed EM scheme will be established. In section $4,$ the almost sure exponential stability for  tamed EM scheme for NSDDEs will be developed.

%In Section $4,$ some auxiliary results to NSDDEs driven by pure jump process will be illustrated, 
%
%An example will be provided in the Section $5$ for demonstration.

\section{Preliminaries}
Throughout this paper, let $(\Omega,\mathcal{F},\PP)$ be a complete probability space with a filtration $\{\mathcal{F}_t\}_{t\geq 0}$ satisfying the usual condition (i.e. it is right continuous and $\mathcal{F}_{0}$ contains all $\PP$-null sets). Let $|\cdot|$ be the Euclidean norm in $\mathbb{R}^n.$ If $A$ is a matrix, denote $||A|| = \sqrt{\trace(A^TA)}$ as the Hilbert-Schmidt norm. Let $C(\mathbb{R}^n;\mathbb{R}_+)$ denote the family of all continuous function from $\mathbb{R}^n$ to $\mathbb{R}_+.$  Denote $\lfloor k\rfloor$ as the integer part of a real number $k.$
%Let $C^2(\mathbb{R}^n;\mathbb{R}_+)$ denote the family of all $C^2$ functions $V$ from $\mathbb{R}^n$ to $\mathbb{R}_+.$

In this paper, we consider the following $n$-dimensional NSDDE 
\begin{equation}\label{NSDDE}
d[X(t)-D(X(t-\tau))]  = b(X(t),X(t-\tau))dt+\sigma(X(t),X(t-\tau))dw(t),\quad t \geq 0
\end{equation}
with initial data satisfies the following condition: for any $p\ge 2$
\begin{equation}\label{initial}
\{X(\theta):-\tau\leq \theta\leq 0\} = \xi \in L^p_{\mathcal{F}_{0}}([-\tau,0];\mathbb{R}^{n}),
\end{equation}
that is  $\xi $ is an $\mathcal{F}_{0}$-measurable $C([-\tau,0];\mathbb{R}^{n})$-valued random variable and  $E\|\xi\|^p <\infty,$ and $w(t)=(w_1(t),\cdots,w_m(t))^T$ is an $m$-dimensional standard $\mathcal{F}_t$-Wiener process. $D:\mathbb{R}^{n} \to \mathbb{R}^{n},$
$b:\mathbb{R}^{n}\times \mathbb{R}^{n}\to \mathbb{R}^{n},$
$\sigma:\mathbb{R}^{n}\times \mathbb{R}^{n} \to \mathbb{R}^{n\times m}$ are Borel-measurable functions.

To show the existence and uniqueness of solution to \eqref{NSDDE},  we assume that:

\begin{description}
	\item[(A1)]There exists a positive constant $\tilde{K}$ such that
	\begin{equation}\label{WeakCorecivity}
	\langle x-D(y),b(x,y)\rangle \vee |\sigma(x,y)|^{2}\leq \tilde{K}(1+|x|^{2}+|y|^2),
	\end{equation}
	for $\forall$ $x,y \in \mathbb{R}.$\\
	\item[(A2)]  $D(0)=0$ and there exists a constant $\kappa \in (0,1)$ such that
	\begin{eqnarray}\label{NeutralBound}
	|D(x)-D(\bar{x})| \leq \kappa|x-\bar{x}| \mbox{ for all }x,y \in \mathbb{R}^n.
	\end{eqnarray}
	\item[(A3)] For any $R>0,$ there exist two  positive constants $\tilde{K}_{R}$ and $K_{R}$ such that
	\begin{equation}\label{LocalLipschitzII}
	\begin{split}
	& \langle x-D(y)-\bar{x}+D(\bar{y}), b(x,y)-b(\bar{x},\bar{y})\rangle\vee|\sigma(x,y)-\sigma(\bar{x},\bar{y})|^2 \\
	&\leq \tilde{K}_{R} (|x-\bar{x}|^{2}+|y-\bar{y}|^2) ,
	\end{split}
	\end{equation}
	for all $|x|\vee|y|\vee|\bar{x}|\vee|\bar{y}| \leq R,$\\
	and 
	\begin{equation}
	\sup_{|x|\vee |y|\leq R}|b(x,y)| \leq K_{R}.
	\end{equation}
\end{description}

\begin{remark}\label{EUTheorem}
	{\rm Assume that {\bf (A1)}-{\bf (A3)} hold, then NSDDE \eqref{NSDDE} with initial \eqref{initial} admits a unique strong global solution $X(t),$ $t\in[0,T].$ The proof details of such existence and uniqueness as well as the $p$th moment bound results can be found in \cite{Ji1}. }
	%this remark should be rewrite since the A1 here is stronger than the one in Ji2.	
	
\end{remark}

With the existence and uniqueness theorem in hand, we could define the the classic EM scheme to \eqref{NSDDE} with initial data \eqref{initial}. For $T >\tau > 0$, without loss of generality, we assume that the ratio of $\frac{T}{\tau}$ is a rational number, and the step size $h\in (0, 1)$ be fraction of $\tau$ and $T,$  such that there exist two positive integers $M, m$ such that $h=T/M=\tau/m.$
Then the discrete-time EM scheme is given by
\begin{itemize}
	\item For $k = -m,\cdots,0,$
	$$Y_{k} = \xi(kh);$$
	\item For $k = 0, \cdots, M-1,$
	$$Y_{k+1}-D(Y_{k+1-m})= Y_{k}-D(Y_{k-m})+b(Y_{k},Y_{k-m})h+\sigma (Y_{k},Y_{k-m})\Delta w_k,$$  
\end{itemize}
where $\Delta w_k: = w_{(k+1)h}- w_{kh}$ and $Y_{k-m}:=Y((k-m)h).$

According to \cite{Hut}, the classic EM scheme may fail if the drift coefficient grow faster than linear rate. As an alternative explicit numerical scheme, the tamed EM scheme has been introduced by \cite{Jen}. The general idea of taming scheme is to set an bound of drift so that its growth rate could be controlled.

Let
\begin{eqnarray}\label{TamedDrift}
b_{h}(x,y):= \frac{b(x,y)}{1+h^{\alpha}|b(x,y)|},\quad \forall x,y\in \mathbb{R}^{n}, \quad \alpha \in (0,\frac{1}{2}]
\end{eqnarray}
be the tamed drift coefficient. 

For $k = 0, \cdots, M-1,$ the discrete-time tamed EM scheme could be rewritten as
\begin{equation}\label{EM2}
\begin{split}
Y_{k+1}-D(Y_{k+1-m}) &= \xi(0)-D(\xi(-\tau))+\sum_{i=0}^{k}b_h(Y_{i},Y_{i-m})h\\
&+\sum_{i=0}^{k}\sigma (Y_{i},Y_{i-m})\Delta w_k.
\end{split}
\end{equation}
For the simplicity, we could define the corresponding continuous-time tamed EM scheme

\begin{itemize}
	\item For $\theta\in[-\tau,0],$
	\begin{equation}\label{initial1}
	Y(\theta) = \xi(\theta);
	\end{equation}
	
	\item For $t\in[t,T],$ 
	\begin{equation}\label{CEulerscheme}
	\begin{split}
	Y(t)=& D(\bar{Y}(t-\tau))+\xi(0)-D(\xi(-\tau))+\int_{0}^{t}b_{h}(\bar{Y}(s),\bar{Y}(s-\tau))ds\\
	&+\int_{0}^{t}\sigma(\bar{Y}(s),\bar{Y}(s-\tau))dw(s),
	\end{split}
	\end{equation}
\end{itemize}
where $\bar{Y}(s) = Y([\frac{s}{h}]h).$

\begin{remark}\label{EUTheorem1}
	{\rm Under assumptions {\bf (A1)}-{\bf (A3)}, the $p$-th moment bounds and convergence result of tamed EM scheme to equation \eqref{NSDDE} can be found in \cite{Ji2}. }

\end{remark}

\section{Exponential Mean-Square Stability}
%\begin{defn}
%	For $p>0,$ we say that the solution $X(t)$ to \eqref{NSDDE} with initial condition \eqref{initial} is said to be $p$-th moment exponentially stable if 
%	\begin{equation}\label{DefStable1CP}
%		\lim\limits_{t\to \infty} \sup\frac{\log E(|X(t)|^p)}{t} < 0.
%	\end{equation}
%\end{defn}

\begin{defn}
	The solution $Y(t)$ to equation \eqref{CEulerscheme} with inital data \eqref{initial1}  is said to be  exponentially mean-square stable if 
	\begin{equation}\label{DefStable1P}
	\lim\limits_{k\to \infty} \sup\frac{\log E(|Y(t)|^2)}{kh} < 0.
	\end{equation}
\end{defn}

\begin{theorem}
	Let Assumptions {\bf (A1)}, {\bf (A2)} and {\bf (A3)} hold. 
	%	Assume that there exists a function $V\in C^2(\mathbb{R}^n;\mathbb{R}_+)$ and some positive constants $c_1,$ $c_2,$ $p$ such that for any $x\in \mathbb{R}^n,$
	%	\begin{equation}\label{StableCon1}
	%	c_1|x|^p\leq V(x)\leq c_2 |x|^p,
	%	\end{equation} 
	Assume that $\sigma$ satisfies the following condition:
	\begin{equation}\label{StableCon2}
	|\sigma(x,y)|^2\leq \frac{1}{h}(-\lambda_1-\lambda_2|x|^2+\lambda_3|y|^2),
	\end{equation}
	where $\lambda_1>2$ and $\lambda_2>\lambda_3>0.$ Then there exists a constant $\tilde{h}\in (0,1),$ such that for any $0<h<\tilde{h}<1,$ the tamed EM approximate solution is said to be  exponentially mean-square stable.
\end{theorem}

\textbf{Proof}: 
%Firstly, we show the stability for the discrete version of tamed EM scheme \eqref{EM2}.
For simplicity, denote that $Z_k = Y_k - D(Y_{k-m}),$  then
\begin{equation}\label{eq1}
\begin{split}
|Z_{k+1}|^2 &= |Z_{k}+b_h(Y_k,Y_{k-m})h+\sigma(Y_k,Y_{k-m})\Delta w_k|^2\\
& = |Z_k|^2+|b_h(Y_k,Y_{k-m})h|^2 +|\sigma(Y_k,Y_{k-m})\Delta w_k|^2\\
&+2\langle Z_k,b_h(Y_k,Y_{k-m})h\rangle+2\langle Z_k,\sigma(Y_k,Y_{k-m})\Delta w_k\rangle\\
&+2\langle b_h(Y_k,Y_{k-m})h, \sigma(Y_k,Y_{k-m})\Delta w_k\rangle.\\
%&\leq 2|Z_k|^2+2|b_h(Y_k,Y_{k-m})h|^2 +|\sigma(Y_k,Y_{k-m})\Delta W_k|^2\\
%&+2\langle Z_k,\sigma(Y_k,Y_{k-m})\Delta w_k\rangle+2\langle b_h(Y_k,Y_{k-m})h, \sigma(Y_k,Y_{k-m})\Delta w_k\rangle,\\
\end{split}
\end{equation}
%where the following inequality has been applied
%\begin{equation*}
%2\langle Z_k,b_h(Y_k,Y_{k-m})h\rangle \leq |Z_k|^2 +|b_h(Y_k,Y_{k-m})h|^2.
%\end{equation*}

By an application of condition (A1), we have 
\begin{equation*}
\langle Z_k,b_h(Y_k,Y_{k-m})h\rangle \leq \tilde{K} (1+|Y_k|^2+|Y_{k-m}|^2).
\end{equation*}

According to the definition of tamed scheme, we have  
\begin{equation*}
|b_h(Y_k,Y_{k-m})h|^2 \leq h^{2-2\alpha}.
\end{equation*}

By letting $|\sigma(x,y)|^2\leq h(-\lambda_1-\lambda_2|x|^2+\lambda_3|y|^2)$ and rewrite 
\begin{equation*}
\begin{split}
|\sigma(Y_k,Y_{k-m})\Delta w_k|^2 &= |\sigma(Y_k,Y_{k-m})|^2(|\Delta w_k|^2-h+h)\\
&\leq |\sigma(Y_k,Y_{k-m})|^2h+|\sigma(Y_k,Y_{k-m})|^2(|\Delta w_k|^2-h)\\
&\leq (-\lambda_1-\lambda_2|Y_k|^2+\lambda_3|Y_{k-m}|^2)h+|\sigma(Y_k,Y_{k-m})|^2(|\Delta w_k|^2-h)
\end{split}
\end{equation*}

the \eqref{eq1} becomes
\begin{equation}\label{eq2}
\begin{split}
|Z_{k+1}|^2&\leq |Z_k|^2+h^{2-2\alpha}+(-\lambda_1-\lambda_2|Y_k|^2+\lambda_3|Y_{k-m}|^2)h\\
&+2\tilde{K}(1+|Y_k|^2+|Y_{k-m}|^2)h+2\langle Z_k,\sigma(Y_k,Y_{k-m})\Delta w_k\rangle\\
&+2\langle b_h(Y_k,Y_{k-m})h, \sigma(Y_k,Y_{k-m})\Delta w_k\rangle+|\sigma(Y_k,Y_{k-m})|^2(|\Delta w_k|^2-h)\\
&=|Z_k|^2+(h^{1-2\alpha}-\lambda_1+2\tilde{K})h+(-\lambda_2+2\tilde{K})|Y_k|^2h\\
&+(\lambda_3+2\tilde{K})|Y_{k-m}|^2h+M^k_1+M^k_2+M^k_3,
\end{split}
\end{equation}
where 
\begin{equation*}
\begin{split}
&M^k_1: = 2\langle Z_k,\sigma(Y_k,Y_{k-m})\Delta w_k\rangle;\\
&M^k_2: = 2\langle b_h(Y_k,Y_{k-m})h, \sigma(Y_k,Y_{k-m})\Delta w_k\rangle;\\
&M^k_3: = |\sigma(Y_k,Y_{k-m})|^2(|\Delta w_k|^2-h).
\end{split}
\end{equation*}

%Note that here, we have an condition on $\lambda_1$, h^{1-2\alpha}-\lambda_1+2\tilde{K}<0

For some positive constant $C>1,$  multiplying $C^{(k+1)h}$ to both sides of \eqref{eq2}, we have
\begin{equation*}
\begin{split}
C^{(k+1)h}|Z_{k+1}|^2&\leq C^{(k+1)h}|Z_k|^2+C^{(k+1)h}(-\lambda_2+2\tilde{K})|Y_k|^2h\\
&+C^{(k+1)h}(\lambda_3+2\tilde{K})|Y_{k-m}|^2h+C^{(k+1)h}(M^k_1+M^k_2+M^k_3),
\end{split}
\end{equation*} 
which is equivalent to 
\begin{equation*}
\begin{split}
C^{(k+1)h}|Z_{k+1}|^2-C^{kh}|Z_k|^2&\leq (C^{(k+1)h}-C^{kh})|Z_k|^2+C^{(k+1)h}(-\lambda_2+2\tilde{K})|Y_k|^2h\\
&+C^{(k+1)h}(\lambda_3+2\tilde{K})|Y_{k-m}|^2h+C^{(k+1)h}(M^k_1+M^k_2+M^k_3).
\end{split}
\end{equation*}

By taking the summation of both sides, we have
\begin{equation*}
\begin{split}
C^{(k+1)h}|Z_{k+1}|^2-C^{0h}|Z_k|^2&\leq \sum_{i=0}^{k}[(C^{(i+1)h}-C^{ih})|Z_i|^2+C^{(i+1)h}(-\lambda_2+2\tilde{K})|Y_i|^2h\\
&+C^{(i+1)h}(\lambda_3+2\tilde{K})|Y_{i-m}|^2h+C^{(i+1)h}(M^i_1+M^i_2+M^i_3)],
\end{split}
\end{equation*}
which implies that
\begin{equation*}
\begin{split}
C^{(k+1)h}|Z_{k+1}|^2&\leq |Z_0|^2+ \sum_{i=0}^{k}[(C^{(i+1)h}-C^{ih})|Z_i|^2+C^{(i+1)h}(-\lambda_2+2\tilde{K})|Y_i|^2h\\
&+C^{(i+1)h}(\lambda_3+2\tilde{K})|Y_{i-m}|^2h+C^{(i+1)h}(M^i_1+M^i_2+M^i_3)],
\end{split}
\end{equation*}

It is not difficult to verify that, for any $i=0,\dots,k,$ $Z_i$ is $\mathcal{F}_{ih}$-measurable. Similarly, $b_{h}(Y_i,Y_{i-m})$ and $\sigma(Y_i,Y_{i-m})$ are also $\mathcal{F}_{ih}$-measurable for all  $i=0,\dots,k.$
Hence $\sum_{i=0}^{k}C^{(i+1)h}(M^i_1+M^i_2)$ is a martingale. 

We also have the martingale properties for the terms involving $M^i_3$  
\begin{equation*}
\begin{split}
E[\sum_{i=0}^{k}C^{(i+1)h}M^i_3\big|F_{kh}] &= E[\sum_{i=0}^{k}C^{ih}M^i_3]+E[C^{(k+1)h}|\sigma(Y_k,Y_{k-m})|^2(|\Delta w_k|^2-h)\big|F_{kh}]\\
& = E[\sum_{i=0}^{k}C^{ih}M^i_3]+C^{(k+1)h}|\sigma(Y_k,Y_{k-m})|^2E[(|\Delta w_k|^2-h)\big|F_{kh}]\\
& = E[\sum_{i=0}^{k}C^{ih}M^i_3],
\end{split}
\end{equation*}
in the last line, we use the fact that 
\begin{equation*}
E[(|\Delta w_k|^2-h)\big|F_{kh}] = 0.
\end{equation*}

Therefore, we have the fact that $\sum_{i=0}^{k}C^{(i+1)h}(M^i_1+M^i_2+M^i_3)$ is a martingale with $M^0_1+M^0_2+M^0_3 = 0.$

According to Lemma 6.4.1 in \cite{Mao}, we have
\begin{equation*}
|Y_i-D(Y_{i-m})|^2\leq (1+\xi)(|Y_i|^2+\frac{|D(Y_{i-m})|^2}{\xi}).
\end{equation*}
By taking $\xi = \kappa^2,$ together with assumption (A2), we have 
\begin{equation*}
|Y_i-D(Y_{i-m})|^2 \leq (1+\kappa^2)(|Y_i|^2+|Y_{i-m}|^2).
\end{equation*}
Then, we could substitute it back to \eqref{eq2} that 
\begin{equation*}
\begin{split}
C^{(k+1)h}|Z_{k+1}|^2&\leq |Z_0|^2+ \sum_{i=0}^{k}\big\{[(C^{(i+1)h}-C^{ih})(1+\kappa^2)+C^{(i+1)h}(-\lambda_2+2\tilde{K})h]|Y_i|^2\\
&+[(C^{(i+1)h}-C^{ih})(1+\kappa^2)+C^{(i+1)h}(\lambda_3+2\tilde{K})h]|Y_{i-m}|^2\\&+C^{(i+1)h}(M^i_1+M^i_2+M^i_3)\big\},
\end{split}
\end{equation*}

%Then our task is to decompose the term of |Y_k-m| into initial data as well as the |Y_k|

Now, we consider the term involving the delay argument, we have 
\begin{equation*}
\begin{split}
&\sum_{i=0}^k[(C^{(i+1)h}-C^{ih})(1+\kappa^2)+C^{(i+1)h}(\lambda_3+2\tilde{K})h]|Y_{i-m}|^2\\
& = \sum_{i=0}^kC^\tau[(C^{(i-m+1)h}-C^{(i-m)h})(1+\kappa^2)+C^{(i+1-m)h}(\lambda_3+2\tilde{K})h]|Y_{i-m}|^2\\
& = C^{\tau}\sum_{i=-m}^{k-m}[(C^{(i+1)h}-C^{ih})(1+\kappa^2)+C^{(i+1)h}(\lambda_3+2\tilde{K})h]|Y_{i}|^2\\
&=C^{\tau}\sum_{i=-m}^{-1}[(C^{(i+1)h}-C^{ih})(1+\kappa^2)+C^{(i+1)h}(\lambda_3+2\tilde{K})h]|\xi(ih)|^2\\
&+C^{\tau}\sum_{i=0}^{m-k}[(C^{(i+1)h}-C^{ih})(1+\kappa^2)+C^{(i+1)h}(\lambda_3+2\tilde{K})h]|Y_{i}|^2\\
&\leq C^{\tau}\sum_{i=-m}^{-1}[(C^{(i+1)h}-C^{ih})(1+\kappa^2)+C^{(i+1)h}(\lambda_3+2\tilde{K})h]|\xi(ih)|^2\\
&+C^{\tau}\sum_{i=0}^{m}[(C^{(i+1)h}-C^{ih})(1+\kappa^2)+C^{(i+1)h}(\lambda_3+2\tilde{K})h]|Y_{i}|^2\\
\end{split}
\end{equation*}

Substitute it back to previous equations, which yields
\begin{equation*}
\begin{split}
&C^{(k+1)h}|Z_{k+1}|^2\leq |Z_0|^2+C^\tau\sum_{i=-m}^{-1}[(C^{(i+1)h}-C^{ih})(1+\kappa^2)+C^{(i+1)h}(\lambda_3+2\tilde{K})h]|\xi(ih)|^2\\ &+\sum_{i=0}^{k}\big\{[(1+C^\tau)(C^{(i+1)h}-C^{ih})(1+\kappa^2)+C^{(i+1)h}(-\lambda_2+C^\tau\lambda_3+2\tilde{K}(1+C^\tau))h]|Y_{i}|^2\\
&+C^{(i+1)h}(M^i_1+M^i_2+M^i_3)\big\}
\end{split}
\end{equation*}

Now, let $f(x) = (1+x^\tau)(x-1)(1+\kappa^2)+(-\lambda_2+x^\tau\lambda_3+2\tilde{K}(1+x^\tau))h$

Obviously, when $x=1,$ $f(1)<0.$ Also, it is clear that $f(x)$ is a continuous function, which is increasing w.r.t $x$. Hence there exists a constant $\bar{C}>1,$ such that $f(\bar{C})=0.$ Therefore, we could choose $C \in (1,\bar{C}),$ such that 
%Here, we need the extra condition that -\lambda_2+\lambda_3+4\tilde{K}<0

\begin{equation}\label{eq3}
\begin{split}
C^{(k+1)h}|Z_{k+1}|^2&\leq |Z_0|^2+C^\tau\sum_{i=-m}^{-1}[(C^{(i+1)h}-C^{ih})(1+\kappa^2)+C^{(i+1)h}(\lambda_3+2\tilde{K})h]|\xi(ih)|^2\\ 
&+C^{(i+1)h}(M^i_1+M^i_2+M^i_3).
\end{split}
\end{equation}

Applying Lemma 6.4.1 from \cite{Mao} again, we have for any $\epsilon>0,$ by choosing $\epsilon = \kappa^2,$
\begin{equation}\label{eq4}
\begin{split}
|Y_{k+1}|^2 &=|Y_{k+1}-D(Y_{k-m+1})+ D(Y_{k-m+1})|\\
&\leq (1+\epsilon)(|Z_{k+1}|^2+\frac{|D(Y_{k-m+1})|^2}{\epsilon})\\
&=(1+\kappa^2)(|Z_{k+1}|^2+|Y_{k-m+1}|^2).
\end{split}
\end{equation}

%Further more, we could have
%\begin{equation*}
%\begin{split}
%|Y_{k+1}|^2&\leq (1+\kappa^2)(|Z_{k+1}|^2+|Y_{k-m+1}|^2)\\
%\end{split}
%\end{equation*}
%The following result could be written as a lemma up front.
%Recall the following fact:  
%\begin{equation*}
%a_k\leq b+qa_{k+1-m}\leq b\sum_{i=0}^{[k/m]}q^i+q^{[(k+1/m)]}|\xi[(k+1-([k+1/m]+1)m)h]|
%\end{equation*}

Taking the expectation of \eqref{eq3}, together with \eqref{eq4}, we have
\begin{equation}\label{eq5}
\begin{split}
&C^{(k+1)h}E|Y_{k+1}|^2\leq (1+\kappa^2)\{E|Z_0|^2\\
&+C^\tau\sum_{i=-m}^{-1}[(C^{(i+1)h}-C^{ih})(1+\kappa^2)+C^{(i+1)h}(\lambda_3+2\tilde{K})h]E|\xi(ih)|^2\}\\
&+(1+\kappa^2)C^\tau C^{(k+1-m)h}E|Y_{k+1-m}|^2.
\end{split}
\end{equation}
Define that

$a_{k+1} := C^{(k+1)h}E|Y_{k+1}|^2,$ $$c_1 :=(1+\kappa^2)\{E|Z_0|^2\\
+C^\tau\sum_{i=-m}^{-1}[(C^{(i+1)h}-C^{ih})(1+\kappa^2)+C^{(i+1)h}(\lambda_3+2\tilde{K})h]E|\xi(ih)|^2\},$$
$c_2:= (1+\kappa^2)C^\tau,$ then equation \eqref{eq5} could be written in the following recursive form:
\begin{equation*} 
\begin{split}
a_{k+1}&\leq c_1+c_2a_{k+1-m}\\
&\leq c_1+c_1c_2+c_2^2a_{k+1-2m}\\
&\leq \dots \leq c_1\sum_{i=0}^{\lfloor k+1/m\rfloor}c_2^i+c_2^{\lfloor k+1/m\rfloor+1}a_{k+1-(\lfloor k+1/m\rfloor+1)m}\\
&= c_1\sum_{i=0}^{\lfloor k+1/m\rfloor}c_2^i+c_2^{\lfloor k+1/m\rfloor+1}E|\xi[(k+1-(\lfloor k+1/m\rfloor+1)m)h]|^2.
\end{split}
\end{equation*}

%here we need to ddfine the notation of the integer part of a real number
Then, there exists a constant $\bar{K}>0,$ which is independent of $k,$ such that
\begin{equation*}
C^{(k+1)h}E|Y_{k+1}|^2 \leq \bar{K}<\infty,
\end{equation*}
which yields that 
\begin{equation*}
\lim\limits_{k\to\infty}\sup\frac{\log E|Y_{k+1}|^2}{(k+1)h}\leq -\log C\log \bar{K}<0.\quad \Box
\end{equation*}

\section{Almost Sure Stability}

\begin{defn}
	The solution $Y(t)$ to equation \eqref{CEulerscheme} with initial data \eqref{initial1} is said to be almost surely exponentially stable if 
	\begin{equation}\label{DefStable1}
	\lim\limits_{k\to \infty} \sup\frac{\log |Y(t)| }{kh}\leq 0\quad a.s.
	\end{equation}
\end{defn}
\begin{coro}
	Under same conditions given in theorem 3.1, the tamed EM scheme is also almost surely stable, i.e. the exponential mean-square stability implies the almost sure stability.
\end{coro}

\textbf{Proof}:
In order to reproduce the almost sure stability, we need to use some result in the proof of exponential mean-square stability. Recall the fact that  there exists a constant $\bar{K}>0,$ which is independent of $k,$ such that
\begin{equation*}
C^{(k+1)h}E|Y_{k+1}|^2 \leq \bar{K}<\infty.
\end{equation*}
By the virtue of Markov inequality, for any $\delta\in(0,h),$ we have 

\begin{equation*}
P(|Y_{k+1}|^2>C^{-(k+1)(h-\delta)})\leq \frac{\bar{K}C^{-(k+1)h}}{C^{-(k+1)(h-\delta)}} = \bar{K}C^{-(k+1)\delta}
, \quad \forall k\geq 1.
\end{equation*}
In the view of Borel-Cantelli lemma yields: for almost all $\omega \in \Omega,$
\begin{equation}\label{eq6}
|Y_{k+1}|^2\leq \bar{K}C^{-(k+1)(h-\delta)}
\end{equation}
holds for all but finitely many $k.$

Furthermore, there exists a $k=k_0(\omega),$ for all $\omega \in \Omega$ excluding a $P$-null set, \eqref{eq6} holds whenever $k\geq k_0.$

Consequently, for almost all $\omega \in \Omega,$
\begin{equation*}
\frac{1}{(k+1)h}\log|Y_{k+1}|= \frac{1}{2(k+1)h}\log|Y_{k+1}|^2\leq -\frac{\log \bar{K}}{2(k+1)h}+\frac{-\log C^{(k+1)(h-\delta)}}{2(k+1)h},
\end{equation*}
whenever $k\geq k_0.$ 
Then 
\begin{equation*}
\lim\limits_{k\to \infty}\sup\frac{1}{(k+1)h}\log|Y_{k+1}|\leq -\frac{\log C}{2}\quad a.s.
\end{equation*}
as $\delta \downarrow 0.$ $\Box$

\begin{remark}
	In general, the exponential mean-square stability doesn't imply the almost sure stability. However, in this case, we notice that once the exponential mean-square stability for tamed EM scheme has been established, the almost sure stability can be reproduced.
\end{remark}

%
%\newpage

\end{document}